\documentclass[12pt]{amsart}
\usepackage{amsmath,amssymb}

\newcommand{\TG}[1]{\mathbb{T}{#1}}
\newcommand{\DG}[1]{\mathbb{D}{#1}}

\newcommand{\E}{\ensuremath{\mathcal E}}
\newcommand{\F}{\ensuremath{\mathcal F}}
\newcommand{\GL}{\textup{GL}}
\newcommand{\GLK}{\textup{GL}_K}

\newcommand{\N}{\ensuremath{\mathbb N}}

\renewcommand{\P}{\ensuremath{\mathcal P}}

\newcommand{\R}{\ensuremath{\mathbb R}}
\newcommand{\RP}{\ensuremath{\mathbb RP}}
\renewcommand{\S}{\ensuremath{\mathbb S}}
\newcommand{\V}{\ensuremath{\mathcal V}}
\newcommand{\W}{\ensuremath{\mathcal W}}

\newcommand{\into}{\hookrightarrow}
\newcommand{\indlim}{\lim\limits_{\longrightarrow}}

\newcommand{\dnc}{{\operatorname{DNC}}}

\theoremstyle{plain}
\newtheorem{thm}{Theorem}[section]

\newtheorem{cor}[thm]{Corollary}
\newtheorem{dfn}[thm]{Definition}
\newtheorem{deflem}[thm]{Definition-Lemma}
\newtheorem{lem}[thm]{Lemma}

\newtheorem{prop}[thm]{Proposition}

\newtheorem{thm*}{Theorem}

  {\addvspace{\bigskipamount}\noindent{\bf Example\ }}%
  {}

\theoremstyle{definition}

\newtheorem{note}[thm]{Note}

\newtheorem{exs}[thm]{Examples}

\begin{document}

\title[Deformation Spaces and Filtrations of Banach Manifolds]
{On Deformation Spaces, Tangent Groupoids and Generalized Filtrations of Banach and Fredholm Manifolds}
\author[Haj Saeedi Sadegh and Jody Trout]{Ahmad Reza Haj Saeedi Sadegh and Jody Trout}
\date{\today}
\email{Ahmadreza.Hajsaeedisadegh@dartmouth.edu, jody.trout@dartmouth.edu}       
\address{6188 Kemeny Hall\\ 
         Dartmouth College\\
         Hanover, NH 03755}


\begin{abstract}
We extend the deformation to the normal cone and tangent groupoid constructions from finite-dimensional manifolds to infinite-dimensional Banach and Fredholm manifolds. Next, we generalize the concept of Fredholm filtrations to get a more flexible and functorial theory. In particular, we show that if $M$ is a Banach (or Fredholm) manifold with generalized filtration $\F = \{M_n\}_1^\infty$ by finite-dimensional submanifolds, then there are induced generalized filtrations $T\F = \{TM_n\}_1^\infty$ of the tangent bundle $TM$ and $\TG{\F} = \{\TG{M_n}\}_1^\infty$ of the tangent groupoid $\TG{M}$, which is not possible in the classical theory.

\end{abstract}

\maketitle

\section{Introduction}

The outline of this paper is as follows. In section \ref{s:deformationspacefredholmtransversality}, we extend some well-known constructions from finite-dimensional smooth manifolds to infinite-dimensional Banach and Fredholm manifolds, such as the deformation to the normal cone $\dnc(M, M_0)$ and Connes' tangent groupoid $\TG{M}$, which have many applications, for example, in the index theory of elliptic operators \cite{Higson2010,HigsonSadegh2018} the study of Lie groupoids \cite{DebordSkandalis21}, quantization theory \cite{QuantCotBdl}, and Connes' noncommutative geometry \cite{Co94}. We also prove some transversality results that are needed for generalized filtrations.

In section \ref{s:generalizedfintration}, we develop the theory of generalized filtrations of Banach and Fredholm manifolds. In the classical theory \cite{Mkhr2, Mkhr, Mkhr3, Mkhr4}, the topology of infinite-dimensional Fredholm manifolds was studied by filtrating the manifold by an increasing sequence $\{M_n\}_{n=k}^\infty$ of finite-dimensional submanifolds whose dimensions $\dim(M_n) = n$ increase without bound and whose union is dense in $M$. The most important condition for algebraic topology is that if $M_\infty = \bigcup_n M_n$ is given the inductive limit topology, then the inclusion $M_\infty \into M$ is a homotopy equivalence. However, a problem with requiring that $\dim(M_n) = n$ is that such a filtration of $M$ does not naturally induce corresponding filtrations of the tangent bundle $TM$ or the pair groupoid $M \times M$ or the tangent groupoid $\TG{M}$ of $M$. A. One would like to have that if $\F = \{M_n\}_1^\infty$ is a filtration  of $M$, then it follows that $T\F = \{TM_n\}_1^\infty$ is a filtration of the tangent bundle $TM$,  $\F_2 = \{M_n \times M_n\}_1^\infty$ is a filtration of the pair groupoid $M \times M$ and $\TG{\F} = \{\TG{M_n}\}$ is a filtration of the tangent groupoid $\TG{M}$. By relaxing the dimension constraint and the density constraint (see Definition \ref{def:Fredfilt}), we develop a theory of generalized filtrations that is more flexible and has these functorial properties, which we illustrate with several examples. 

The first author is grateful to the Max Planck Institute for Mathematics in Bonn for their financial support and for providing office space during June and July of 2025, where part of this work was completed. The second author would like to thank Dorin Dumitra\c scu, Nigel Higson and the referee for helpful suggestions.

\section{Deformation Spaces and Fredholm Manifolds}\label{s:deformationspacefredholmtransversality}

We assume throughout that $\E$ denotes a real separable infinite-dimensional Banach space. All manifolds modeled on $\E$ are assumed to be $C^\infty$, separable, metrizable (hence, paracompact), connected and admit partitions of unity. Partitions of unity will exist if, for example, $\E$ has an equivalent norm that is $C^\infty$ away from the origin. All maps are assumed to be $C^\infty$. For the basic theory of Banach manifolds, see the books by Lang \cite{Lang} and Abraham, Marsden and Ratiu \cite{AMR88}. 

By a pair of Banach manifolds $(M,M_0)$ we mean a Banach manifold $M$ with a closed embedded submanifold $M_0\into M$ modeled on a closed subspace $\E_0 \subseteq \E$ with closed complement. 

In this section, we introduce the deformation to the normal cone $\dnc(M,M_0)$ of pairs of Banach manifolds $(M,M_0)$ that is a Banach manifold with a fibration over the real numbers. Indeed, the $\dnc$ construction is a covariant functor from the category of pairs of smooth Banach manifolds to the category of smooth Banach manifolds that are fibered over $\mathbb{R}$. See \cite{DebordSkandalis21, Higson2010} for an introduction to the deformation to the normal cone in finite-dimensional differential topology.

When the pair consists of Fredholm manifolds, we show that the deformation to the normal cone becomes a Fredholm manifold in a canonical way. As a consequence, we see that the tangent groupoid $\mathbb{T}M$ of a Fredholm manifold $M$ is itself a Fredholm manifold. Finally, we will show that the $\dnc$ functor and, hence, the tangent groupoid functor, respect transversality, which we need to discuss generalized Fredholm filtrations of the tangent groupoid in the next section.

\subsection{Deformation Spaces}
Let $(M, M_0)$ be a pair of smooth Banach manifolds.
The {\bf deformation to the normal cone} $\dnc(M,M_0)$ is, as a set, given by 
\[\nu(M,M_0)\times\{0\} \sqcup \bigsqcup_{\lambda\neq 0 \in \R} M\times\{\lambda\},\]
where $\nu(M,M_0)$ is the normal bundle  $$\pi:\nu(M,M_0):=TM|_{M_0}/TM_0\to M_0$$
of the embedding $M_0\into M$. Note that $\dnc(M,M_0)$ is fibered over the reals $\R$.

The topology (and the manifold structure) of $\dnc(M,M_0)$ is determined by declaring the following maps to be continuous:

\begin{itemize}
    \item The canonical map $p:\dnc(M,M_0)\to M\times\mathbb{R}$ given by
    \[\begin{cases}
        (X,m_0,0) \mapsto (m_0,0)\in M\times \{0\}, & (X,m)\in \nu(M,M_0)|_{m_0}\\
        (m,\lambda)\mapsto (m,\lambda),& \lambda\neq 0.
    \end{cases}\]


    \item For any tubular neighborhood embedding $\phi: \mathbb{U}\to U$ from a neighborhood $\mathbb{U}\subset\nu(M,M_0)$ of the zero section to a neighborhood $U$ of $M_0$ in $M$, the induced map $\tilde{\phi}:\tilde{\mathbb{U}}\to \dnc(M,M_0)$
    \[\begin{cases}
        (X,m_0,0)\mapsto (X,m_0,0), & (X,m)\in \nu(M,M_0)|_{m_0}\\
        (X,m,\lambda)\mapsto (\phi(\lambda X,m),\lambda), &\lambda\neq 0,
    \end{cases}\]
    where $\tilde{\mathbb{U}}:=\{(X,m,\lambda)\in \nu(M,M_0)\times \mathbb{R}: (\lambda X,m)\in\mathbb{U}\}.$

\end{itemize}
Note that the map $p$ from the first bullet point will be a homeomorphism away from the zero fiber, i.e., $p:\dnc(M,M_0)\backslash \nu(M,M_0)\times\{0\} \cong M\times (\R\backslash\{0\}).$ The second bullet point identifies neighborhoods of the zero fiber of $\dnc(M,M_0)$ with the open subsets $\tilde{\mathbb{U}}$ of $\nu(M,M_0)\times\mathbb{R}$. This way $\dnc(M,M_0)$ canonically becomes a smooth manifold.

    In \cite{HigsonSadegh2018}, the topology of $\dnc$ is given in a more canonical way using the notion of $C^{\infty}$-rings, which is similar to the definition given in algebraic geometry  (see \cite{Fulton1998}); however, we will not take on this approach. 

    Indeed, $\dnc$ is a covariant functor from the category of pairs of manifolds to the category of smooth manifolds with a fibration over $\mathbb{R}$. For a morphism of pairs $(f,f_0):(M,M_0)\to (N,N_0)$ we obtain a morphism 
    
    \[\dnc(f):\dnc(M,M_0)\to \dnc(N,N_0)\]
    \[\begin{cases}
        (X,m_0,0)\mapsto ([Df(X)],f_0(m_0),0), & (X,m_0)\in \nu(M,M_0)|_{m_0}\\
        (m,\lambda)\mapsto (f(m),\lambda), & \lambda\neq 0.
    \end{cases}\]

    We are interested in a special case of the deformation to the normal cone given by the tangent groupoid.

    \begin{dfn}
        Let $M$ be a Banach manifold. The tangent groupoid of $M$, denoted by $\TG M$, is the deformation space $\dnc(M\times M, M)$ associated with the diagonal embedding $M\into M\times M$. So, as a set, we have
        \[\TG M=TM\times\{0\} \sqcup\bigsqcup_{\lambda\neq 0 \in \R} M\times M\times\{\lambda\},\]
        where we used the identification $\nu({M\times M},M)\cong TM$ given by
        \[\nu(M\times M,M)=T(M\times M)|_{M}/TM\to TM\]
        \[[v_1,v_2]\mapsto v_1-v_2.\]
        The tangent groupoid can be considered a covariant functor from the category of smooth Banach manifolds to the category of smooth Banach manifolds fibered over $\mathbb{R}$ (also to the category of Banach-Lie groupoids). For a morphism $f:M\to N$ we write the corresponding morphism on the $\mathbb{T}$ functor by $\mathbb{D}f:\mathbb{T}M\to\mathbb{T}N.$
    \end{dfn}

The following lemma is quite elementary to check.

\begin{lem}\label{l:dnccartesianproduct}
    Let $(M,M_0)$ and $(N,N_0)$ be pairs of Banach manifolds. Then we have a canonical diffeomorphism
    \[\dnc(M\times N,M_0\times N_0)\cong \dnc(M,M_0)\times_{\mathbb{R}}\dnc(N,N_0).\]
    In particular, for Banach manifolds $M,N$, we have a canonical diffeomorphism
    \[\TG (M\times N)\cong \TG M\times_{\R}\TG N.\]
\end{lem}

\begin{lem}\label{l:dncofbanachspace}
    Let $\E$ be a Banach space with $\E_0\subseteq \mathcal{E}$ a closed complemented subspace. Then $\dnc(\mathcal{E},\E_0)$ is isomorphic to $\mathcal{E}\times \R$ as Banach manifolds over $\R$.
    In particular, $\TG{}\mathcal{E}$ is diffeomorphic to $\mathcal{E}\times \mathcal{E}\times\R$ as Banach manifolds fibered over $\R$.
\end{lem}

\begin{proof}
    Let $V$ be a complement for $\E_0$ in $\mathcal{E}$; so for every $e\in \mathcal{E}$ there exist unique vectors $e'\in \E_0$ and $v\in V$ such that $e=e'+v$. Now we have the isomorphism 
    \[\dnc(\mathcal{E},\E_0)\to \mathcal{E}\times \mathbb{R}\]
    \[\begin{cases}
        (e'+v,t)\mapsto (e'+\frac{v}{t},t) & t\neq 0\\
        (e',v,0)\mapsto (e'+v,0)& t=0.
    \end{cases}\]
\end{proof}

\begin{prop}\label{prop:trivialbdl}
    Let $M\times \mathcal{E}\to M$ be a trivial bundle, then $\TG (M\times \mathcal{E})$ is isomorphic to $\TG M\times (\mathcal{E}\times \mathcal{E})$ as a vector bundle over $\TG M$.
\end{prop}

\begin{proof}
    By Lemma \ref{l:dnccartesianproduct} and Lemma \ref{l:dncofbanachspace}, we have
    \[\TG (M\times\mathcal{E})\cong \TG M\times_{\R} \TG{}\mathcal{E}\cong \TG M\times_{\R} (\mathcal{E}\times\mathcal{E}\times\R)\cong \TG M\times (\mathcal{E}\times\mathcal{E}).\]
    These diffeomorphisms are all isomorphisms of vector bundles over $\TG M$.
\end{proof}

\subsection{DNC-functor and Fredholm structures}
We denote by $\GL(\E)$ the group of bounded invertible operators on $\E$ and by $K(\E)$ the space of all compact operators on $\E$. We let  $\GLK(\E)$ denote the closed normal subgroup of $\GL(\E)$ given by 
\[\GLK(\E):=\{ T = I + K  \in \GL(\E) : K \in K(\E)\} = (I + K(\E)) \cap \GL(\E),\]
which is called the {\it Fredholm structure group}.

Let $M$ be a smooth Banach manifold modeled on $\E$. We say $M$ carries a Fredholm structure, or simply $M$ is a \emph{Fredholm manifold}, if there exists an atlas compatible with the given smooth structure such that the differentials of the transition functions belong to the Fredholm structure group $\GLK(\E)$. This means that for every two charts $\phi_i:U_i\to \E$ of the atlas with $U_i \cap U_j \neq \emptyset$, the transition function
\[\phi_2\circ\phi_1^{-1}:\phi_1(U_1\cap U_2)\to \phi_2(U_1\cap U_2),\]
has differential
\[D(\phi_2\circ\phi_1^{-1})(x)\in\GLK(\E)\]
for all $x\in \phi_1(U_1\cap U_2).$ For a more detailed discussion on the Fredholm structure group and Fredholm manifolds, along with examples, see \cite{ElwTr, EeElw1, Mkhr, Mkhr3, DT06, Ksch}. If $\GL(\E)$ is contractible (e.g., if $\E$ is a Hilbert space), then any Banach manifold modeled on $\E$ can be given a Fredholm structure. However, even for a Hilbert manifold, there can be infinitely many distinct Fredholm structures (see Example 6.2 \cite{Ksch}). 

We want to show for a pair of Fredholm manifolds $(M,M_0)$, the deformation to the normal cone $\dnc(M,M_0)$ carries a canonical Fredholm structure. As a consequence, we deduce for a Fredholm manifold $M$, the tangent groupoid $\mathbb{T}M$ is a Fredholm manifold.

\begin{lem}\label{l:Fredholmindexzero}
    Let $F:\mathcal{E}_1\to \mathcal{E}_2$ and $F':\E'_1\to\E'_2$ be Fredholm operators of index zero and let $P:\mathcal{E}_1\to\mathcal{E}_2'$ be a bounded operator. Suppose $T:\mathcal{E}_1\oplus\E'_1\to \mathcal{E}_2\oplus\E'_2$  has the matrix form
    \[\begin{bmatrix}
        F&0\\
        P&F'
    \end{bmatrix}.\]
    Then $T$ is Fredholm of index zero. 
    \end{lem}

\begin{proof}
    Since $F,F'$ are Fredholm of index zero, there are compact (even finite-rank) operators $K:\mathcal{E}_1\to \mathcal{E}_2$ and $K':\E'_1\to \E'_2$ such that $F+K$ and $F'+K'$ are invertible. Hence $T+(K\oplus K')$ is invertible, which implies $T$ is Fredholm of index zero.
\end{proof}

If $F\in\GLK(\E)$ and $F'\in\GLK(\E')$ then a bounded linear operator $T : \E \oplus \E' \to \E \oplus \E'$ of the matrix form 
 \begin{equation}\label{eq:lowertriangularform}\begin{bmatrix}
        F&0\\
        P&F'
    \end{bmatrix}\end{equation}
is in $\GL(\E \oplus \E')$ (thus, Fredholm of index zero), but $T$ does necessarily belong to the subgroup $\GLK(\E\oplus\E')$. We need the following lemma:

\begin{deflem}\label{dl:lowertriangulardeformationretractdiagonal}
    Let $\widetilde\GL_K(\E\oplus\E')$ denote the group of invertible operators of the lower block triangular form \eqref{eq:lowertriangularform}, where $F\in\GLK(\E)$, $F'\in\GLK(\E')$, and $P:\E\to\E'$ is bounded. Then $\widetilde\GLK(\E\oplus\E')$ deformation retracts to the subgroup $\GLK(\E)\times \GLK(\E')\subset \GLK(\E\oplus\E')$.
\end{deflem}

\begin{proof}
    The proof is explicit. The deformation retraction is given by the map
    \[\Phi_t:\widetilde\GLK(\E\oplus\E')\to \widetilde\GLK(\E\oplus\E'), \ \ \ 0\leq t\leq 1,\]
    \[\begin{bmatrix}
        F&0\\
        P&F'
    \end{bmatrix}\mapsto \begin{bmatrix}
        F&0\\
        (1-t)P&F'
    \end{bmatrix}.\]
\end{proof}

\begin{prop}\label{p:decompositiontangentspacenormalbundle}
    Fix a connection $\nabla$ on the normal bundle; then the tangent bundle of the normal bundle $T\nu(M,M_0)$ has a decomposition 
    \[T\nu(M,M_0)\cong \pi^*TM_0\oplus \pi^*\nu(M,M_0).\]
\end{prop}

This Proposition follows from a more general statement \cite{AMR88,Lang}:
\begin{prop}
    Let $\pi:E\to M_0$ be a vector bundle with a connection $\nabla$, then there is a canonical decomposition
    \[TE\cong \pi^*TM_0\oplus \pi^*E\]
    as vector bundles over $E$.
\end{prop}

Now consider two pairs of Banach manifolds $(M,M_0)$ and $(N,N_0)$. Let $$(f,f_0):(M,M_0)\to (N,N_0)$$ a \emph{map of pairs}, i.e., a smooth map $f:M\to N$ whose restriction to $M_0$ gives a map $f_0:M_0\to N_0$. Such a map induces a morphism of the normal bundle:
\[f_*:\nu(M,M_0)\to\nu(N,N_0)\]
\[f_*([X],m)\mapsto ([Df(X)],f(m)).\]

The following statement follows from straightforward local calculations:

\begin{lem}\label{l:mapnormalbundleindexzero}
    The differential of the induced map on the normal bundle
     \[Df_*:T\nu(M,M_0)\to T\nu(N,N_0)\]
     is of the form 
     \[\begin{bmatrix}
    Df_0&0\\
    P& \pi^*f_*
    \end{bmatrix}\]
    with respect to the decompositions $T\nu(M,M_0)\cong \pi^*TM_0\oplus \pi^*\nu(M,M_0)$ and $T\nu(N,N_0)\cong \pi^*TN_0\oplus \pi^*\nu(N,N_0)$ (as in Proposition \ref{p:decompositiontangentspacenormalbundle}).
    
\end{lem}

\begin{prop}\label{p:normalbundlefredholm}
       If $M$ is a Fredholm manifold with $M_0\subset M$ a Fredholm submanifold, then $\nu(M,M_0)$ is also a Fredholm manifold. In particular, if $M$ is a Fredholm manifold, then $TM$ is also a Fredholm manifold.
\end{prop}
\begin{proof}
    For the first part of the statement, fix a Banach space pair $(\E,\E_0)$ over which the pair $(M,M_0)$ is modeled, and fix a complement $\E'$ of $\E_0$ in $\E$.  Take two local coordinates $\phi_i:U_i\to \mathcal{U}_i$ (i=1,2) from neighborhoods $U_i\subset M$ of a point in $M_0$, to neighborhoods of the origin $\mathcal{U}_i\subset \mathcal{E}$, such that $\psi_i:=\phi_i|_{U_i\cap M_0}$ maps $U_i\cap M_0$ maps diffeomorphically to $\mathcal{U}_i\cap \E_0$. The pair of transition functions 
    \[(f,f_0):=(\phi_2\circ \phi_1^{-1},\psi_2\circ\psi_1^{-1}):(\mathcal{U}_1,\mathcal{U}_1\cap M_0)\to (\mathcal{U}_2,\mathcal{U}_2\cap M_0)\]
    are $\GLK$-maps on $(\E,\E_0)$. This pair induces a map
    \[f_*:\nu(\mathcal{U}_1,\mathcal{U}_1\cap M_0)\to \nu(\mathcal{U}_2,\mathcal{U}_2\cap M_0).\]
    Using Proposition \ref{p:decompositiontangentspacenormalbundle} and Lemma \ref{l:mapnormalbundleindexzero}, with respect to the decomposition $T_x\nu(\mathcal{U}_i,\mathcal{U}_i\cap M_0)\cong \E_0\oplus \E' $, the differential map $Df_*$ is of the form 
    \[\begin{bmatrix}
    Df_0&0\\
    P& \pi^*f_*
    \end{bmatrix}.\]
    Hence, $Df_*\in \widetilde{\GLK}(\E_0\oplus\E')$, i.e., $\nu(M,M_0)$ carries a canonical $\widetilde\GLK(\E_0\oplus\E')$ structure. By the Definition-Lemma \ref{dl:lowertriangulardeformationretractdiagonal}, this structure reduces to a canonical $\GLK(\E_0\oplus\E')$ structure. The second statement follows from the identification of the tangent bundle $TM\cong\nu(M\times M,M).$
    \end{proof}

Now we state the main theorem of this section.

\begin{thm}\label{th:dncfredholm}
If $(M,M_0)$ is a pair of Fredholm manifolds, then $\dnc(M,M_0)$ has a canonical Fredholm structure. In particular, if $M$ is a Fredholm manifold, then the tangent groupoid $\mathbb{T}M$ is canonically a Fredholm manifold.    
\end{thm}
    \begin{proof}
For $t\neq0$, near $(m,t)\in \dnc(M,M_0)$ 
the deformation space looks like $M\times (t-\epsilon,t+\epsilon)$, for $\epsilon$ small. So, the Fredholm structure is clear.

For $t=0$, near $\nu(M,M_0)\times\{0\}$, the deformation space $\dnc(M,M_0)$ is locally diffeomorphic to $\nu(M,M_0)\times [-\epsilon,\epsilon]$; fix a neighborhood $\mathcal{U}\subset \nu(M,M_0)$ and $U\subset M$ of a point $m\in M_0$ with a tubular neighborhood embedding
\[\phi:\mathcal{U}\to U.\] 
For $\epsilon>0$ small enough, we have a map
\[\tilde{\phi}:\mathcal{U}\times[-\epsilon,\epsilon]\to\dnc(M,M_0)\]
\begin{equation}\label{eq:localtrivializationofdnc}
    \begin{cases}
    (X,m,t)\mapsto (\phi(tX),m,t)& t\neq0\\
    (X,m,0)\mapsto (X,m,0)&t=0
\end{cases}
\end{equation}
that is a diffeomorphism onto a neighborhood of $(m,0)\in\dnc(M,M_0)$. From Proposition \ref{p:normalbundlefredholm} it follows that near $\nu(M,M_0)\times\{0\}$, the deformation space $\dnc(M,M_0)$ has a canonical Fredholm structure.
\end{proof}

A smooth map $f : M \to N$ between Banach manifolds is called {\it Fredholm} if the differential $Df_x : T_xM \to T_{f(x)}N$ is a Fredholm operator for all $x \in M$. The Fredholm index of $Df_x$ is independent of $x$ when $M$ is connected, as we assume, which we call the index of the map $f$.

\begin{prop}\label{p:dncmapfredholm}
    Let $(f,f_0):(M,M_0)\to (N,N_0)$ be a smooth map of pairs of Banach manifolds, each Fredholm of index zero. Then the induced maps
    \[f_*:\nu(M,M_0)\to \nu(N,N_0)\]
    and 
    \[\dnc(f):\dnc(M,M_0)\to\dnc(N,N_0)\]
    are Fredholm maps of index zero. In particular, if $f:M\to N$ is a smooth Fredholm map of index zero, then the induced maps
    \[Df:TM\to TN\]
    of the tangent bundles and 
    \[\mathbb{D} f:\mathbb{T} M\to \mathbb{T} N\]
    of the tangent groupoids are also Fredholm maps of index zero.
    
\end{prop}

\begin{proof}
    The statement for $f_*:\nu(M,M_0)\to \nu(N,N_0)$ follows from Lemmas \ref{l:Fredholmindexzero} and \ref{l:mapnormalbundleindexzero}. For the map $\dnc(f)$ the argument is very similar to the proof of Theorem \ref{th:dncfredholm}. 
    The second part of the proposition on the maps $Df$ and $\mathbb{D}f$ follows from the first part.
\end{proof}

\subsection{DNC-functor and Transversality}
In order to construct generalized Fredholm filtrations on the tangent bundle and tangent groupoid of a Fredholm manifold, we need some transversality results in this subsection. First, we show that the normal bundle functor is compatible with transversality. Then, it follows that the same property holds for the $\dnc$-functor. To do so, we need the following lemma.

\begin{lem}\label{l:lowertriangulartransversality}
    Let $T_i:\E_1\to \E'_1$, $i=1,2$, be bounded operators that are transversal to closed complemented subspaces $\mathcal{V}_i\subset\E'_i$, respectively. If $P:\E_1\to\E'_2$ is a bounded operator, then the operator $T:\E_1\oplus\E_2\to\E'_1\oplus\E'_2$ with the matrix form
    \[\begin{bmatrix}
    T_1&0\\
    P& T_2
    \end{bmatrix}\]
    is transversal to $\mathcal{V}_1\oplus\mathcal{V}_2$.
    \end{lem}

    \begin{proof}
        To prove transversality, we need to establish two properties:
        \begin{enumerate}
            \item $T(\E_1\oplus\E'_1)+(\V_1\oplus\V_2)=\E'_1\oplus\E'_2$.

            \item $T^{-1}(\V_1\oplus\V_2)$ is complemented in $\E_1\oplus\E_2$.
        \end{enumerate}

To see the first property, take $(e'_1,e'_2)\in \E'_1\oplus \E'_2$. By transversality property of $T_1,T_2$, we can find $e_i\in \E_i$ and $v_i\in\V_i$, $i=1,2$, such that 
\[e'_i=T_i(e_i)+v_i, \ \ \ i=1,2.\]
Using the transversality property of $T_2$, again, we can find $e\in \E_2$ and $v\in \V_2$ such that
\[P(e_1)=T_2(e)+v.\]
Now we can write
\[(e'_1,e'_2)=T(e_1,e_2-e)+(v_1,v_2-v),\]
from which the first property follows.

To prove the second property, take complementary subspaces $\W_i\subset \E_i$ for $T_i^{-1}(\V_i)$, $i=1,2$. We show $\W_1\oplus\W_2$ is a complementary subspace for $T^{-1}(\V_1\oplus\V_2)$.

For any vector $(e_1,e_2)\in\E_1\oplus\E_2$, by transversality property of $T_1,T_2$, we can can find $w_i\in\W_i$ and $u_i\in T_i^{-1}(\V_i)$, $i=1,2$, such that 
\[e_i=u_i+w_i.\]
Again, using the transversality property of $T_1$, we can find  $u\in \E_2$ and $v\in\V_1$ such that $P(e_1)=T_2(u)+v$. By writing $u=\tilde{v}+\tilde{w}$ where $\tilde{v}\in T_2^{-1}(\V_2)$ and $\tilde{w}\in \W_2$, we have the equality
\[(e_1,e_2)=(u_1,u_2-\tilde{w})+(w_1,w_2+\tilde{w})\]
where $(u_1,u_2-\tilde{w})\in T^{-1}(\V_1\oplus\V_2)$ and  $(w_1,w_2+\tilde{w})\in\W_1\oplus\W_2$. From this  and the fact that $T^{-1}(\V_1\oplus\V_2)\cap (\W_1\oplus\W_2)=\{0\}$ the second property follows.

\end{proof}

\begin{prop}\label{p:trasnversalitynormalmap}
    Let $(f,f_0):(M,M_0)\to (N,N_0)$ be map of pairs of Banach manifolds. Assume $Z\into N$ is a submanifold transverse to $N_0$ and denote $Z_0:=Z\cap N_0$. If $f$ and $f_0$ are transversal to $Z$ and $Z_0$, respectively, then the induced map $f_*:\nu(M,M_0)\to\nu(N,N_0)$ is transversal to $\nu(Z,Z_0)$. Furthermore, $f_*^{-1}(\nu(Z,Z_0))=\nu(f^{-1}(Z),f_0^{-1}(Z_0))$.
\end{prop}

\begin{proof}
    Similar to proof of Lemma \ref{l:mapnormalbundleindexzero}, the differential $Df_*:T\nu(M,M_0)\to T\nu(N,N_0)$ can be represented by a matrix
    \[\begin{bmatrix}
    Df_0&0\\
    P& \pi^*f_*
    \end{bmatrix}.\]
    Now, from Lemma \ref{l:lowertriangulartransversality}, the first part of the statement follows. The equality $$f_*^{-1}(\nu(Z,Z_0))=\nu(f^{-1}(Z),f_0^{-1}(Z_0))$$ can be checked set-theoretically. \end{proof}

\begin{cor}\label{cor:trasveralitytangent}
    If $f:M\to N$ is transversal to $Z\into N$, then $Df:TM\to TN$ is transversal to $TZ \into TN$. Also, $Df^{-1}(TZ)=T(f^{-1}(Z))$.
\end{cor}

\begin{thm}
    Let $(f,f_0):(M,M_0)\to (N,N_0)$ be a map of pairs of Banach manifolds. Assume $Z\into N$ is a submanifold transverse to $N_0$ and denote $Z_0:=Z\cap N_0$. If $f$ and $f_0$ are transverse to $Z$ and $Z_0$, respectively, then the induced map $\dnc(f):\dnc(M,M_0)\to\dnc(N,N_0)$ is transversal to $\dnc(Z,Z_0)$. We also have $$\dnc(f)^{-1}(\dnc(Z,Z_0))=\dnc(f^{-1}(Z),f_0^{-1}(Z_0)).$$
\end{thm}

\begin{proof}
Over the fibers $t\neq0$ of the deformation spaces, the statement is clear as the function $\dnc(f)$ equals $(f,\textup{id}):M\times\mathbb{R}^{*}\to N\times\mathbb{R}^*$ which is transversal to $Z\times\mathbb{R}^*$.    

Near the fiber $t=0$, we may identify the deformation spaces with a product of the normal bundles with an interval as in the proof of Theorem \ref{p:dncmapfredholm}. 
Take neighborhoods $\mathcal{U}_1\subset \nu(M,M_0)$ and $\mathcal{U}_2\subset \nu(N,N_0)$ of points $m\in M_0$ and $f_0(m)\in N_0$, and tubular neighborhood embeddings
\[\phi_i:\mathcal{U}_i\to U_i, \ i = 1, 2, \]
where $U_1\subset M$ and $U_2\subset N$ are open neighborhoods of $m$ and $f_0(m)$ respectively. We then obtain maps 
\[\widetilde\phi_1:\mathcal{U}_1\times[-\epsilon,\epsilon]\to\dnc(M,M_0), \] 
and
\[ \widetilde\phi_2:\mathcal{U}_2\times[-\epsilon,\epsilon]\to\dnc(N,N_0) \]
with formulas similar to \eqref{eq:localtrivializationofdnc}, that are diffeomorphisms onto their images. Then, the composition map
\begin{equation}\label{eq:localformofdncmap}
\widetilde{\phi}_2^{-1}\circ\dnc(f)\circ\widetilde{\phi}_1:\nu(M,M_0)\times[-\epsilon,\epsilon]\to \nu(N,N_0)\times[-\epsilon,\epsilon]    
\end{equation}
is given by the formulas
\[\begin{cases}
    (X,m,t)\mapsto (\frac{1}{t}\phi_2^{-1}\circ f(\phi_1(tX)),f(m),t)& t\neq0\\
    (X,m,0)\mapsto (Df(X),f(m),0)&t=0.
\end{cases}\]
Using Taylor's theorem \cite{Lang}, there exist a smooth function $g:\nu(M,M_0)\times[-\epsilon,\epsilon]\to\nu(N,N_0)$ such that $$\frac{1}{t}\phi_2^{-1}\circ f(\phi_1(tX))=Df(X)+tg(X,t)$$ for $t\neq0$. Hence, we have that $$\widetilde{\phi_2}^{-1}\circ\dnc(f)\circ\widetilde{\phi_1}(X,m,t)=(Df(X)+tg(X,t),f(m),t)$$ for all $t\in[-\epsilon,\epsilon]$. Hence, the differential of the map \eqref{eq:localformofdncmap} is given by a matrix 
\begin{equation*}\label{eq:matrixformofdncmap}
    \begin{bmatrix}
    Df_*& 0\\
    h&\textup{id}
\end{bmatrix}
\end{equation*}
where $h$ denotes a term given in terms of $f,\phi,\psi $ and $t$ and its   derivatives up to order $2$. By Proposition \ref{p:trasnversalitynormalmap} and Lemma \ref{l:lowertriangulartransversality}, it follows that the matrix above is transversal to the preimage $\widetilde\phi_2^{-1}(\dnc(Z,Z_0))$. From this, the theorem follows.
\end{proof}

\begin{cor}\label{cor:transversaltangentgrpoid}
    If $f:M\to N$ is transversal to $Z\into N$, then $\mathbb{D}f:\mathbb{T}M\to \mathbb{T}N$ is transversal to $\mathbb{T}Z \into \mathbb{T}N$, and $\mathbb{D}f^{-1}(\mathbb{T}Z)=\mathbb{T}(f^{-1}Z)$.
\end{cor}


\section{Generalized Filtrations of Banach and Fredholm Manifolds}\label{s:generalizedfintration}

In this section, we develop the theory of generalized filtrations of Banach and Fredholm manifolds and their functorial properties that are not present in the classical theory. We begin by relaxing the dimension requirement in the definition of a flag in a Banach space.

\begin{dfn}\label{def:flag} Let $\Delta=\{\delta(n)\}_{n=1}^\infty$ be a strictly increasing sequence of natural numbers. A {\bf $\Delta$-flag} in a Banach space $\E$ is a sequence of closed subspaces $\{E_n\}_{n=1}^\infty$ of $\E$ such that for all $n \geq 1:$
\begin{itemize}
    \item[a.)] $\dim(E_n) = \delta(n)$.
    \item[b.)] $E_n \subset E_{n+1}$.
    \item[c.)]  $\bigcup_{n=1}^\infty E_n$ is dense in $\E$.
    \item[d.)]$E_n$ has a closed complement $E^{\infty-n}$ in $\E$.
    \item[e.)] $E^{\infty - n} \supseteq E^{\infty - (n+1)}$.
  \end{itemize}
\end{dfn}
A classical flag is an $\N$-flag where $\dim(E_n) = n$ for all $n \geq 1$. Clearly, if $\E$ has a Schauder basis, e.g., if $\E$ is a Euclidean (real Hilbert) space, then $\Delta$-flags exist. Although any $\Delta$-flag can be ``completed'' to an $\N$-flag, the relaxed dimension constraint allows the following constructions, which we will need below.

\begin{lem}\label{lem:flags} If $\{E_n\}_{n=1}^\infty$ is a $\Delta$-flag in the Banach space $\E$, then we have that:
\begin{itemize}
  \item[i.)] Any subsequence $\{E_{n_k}\}_{k=1}^\infty$ is a $\tilde{\Delta}$-flag with $\tilde{\Delta} = \{\delta(n_k)\}_{k=1}^\infty$.
    \item[ii.)] $\{E_n \times E_n\}_{1}^\infty$ is a $2\Delta$-flag in $\E \times \E \cong T\E$.
    \item[iii.)] If $\{F_n\}_1^\infty$ is a $\Delta'$-flag in $\E'$ then $\{E_n \times F_n\}_1^\infty$ is a $(\Delta+\Delta')$-flag in $\E \times \E'$.
    \item[iv.)] $\{E_n \times E_n \times \R\}_1^\infty$ is a $(2\Delta+1)$-flag in $\E \times \E \times \R \cong \TG{\E}$. 
\end{itemize}
\end{lem}
We can now state the main definition of this section.

\begin{dfn}\label{def:Fredfilt}
Let $M$ be a Banach manifold. Let $\Delta=\{\delta(n)\}_{n=1}^\infty$ be a strictly increasing sequence of natural numbers.  A {\bf $\Delta$-filtration} of $M$ consists of a sequence $\F = \{ M_n \}_{n=1}^\infty$ of finite-dimensional closed submanifolds of $M$ satisfying for all $n \geq 1$:
\begin{itemize}
\item[a.)] $\dim(M_n) = \delta(n)$;
\item[b.)] $M_n \subset M_{n+1}$;
\item[c.)] If $M_\infty =  \cup_{n = 1}^\infty M_n$ is given the inductive limit topology, then the  inclusion $i_\infty : M_{\infty} \into M$ is a homotopy equivalence.
\end{itemize}
The $\Delta$-filtration is called {\bf normal} if the following two conditions are also satisfied:
\begin{itemize} 
\item[d.)] The embeddings $M_n \into M_{n+1}$ and
$M_n \into M$ have trivial normal bundles.
\item[e.)] There exist open tubular neighborhoods $V_n$ of $M_n$ in $M$ and open sets $U_n \subseteq V_n$ such that $U_n \subseteq U_{n+1}$ and $M = \cup_{n=1}^\infty U_n$.
\end{itemize}
A $\Delta$-filtration is called {\bf dense} if $M_\infty$ is dense in $M$ and is called {\bf Fredholm} if there is a Fredholm map of index zero $f : M \to \E'$ to some Banach space $\E'$ that is transverse to a $\Delta$-flag $\{E_n\}_{n=1}^\infty$ of $\E'$ for which $M_n = f^{-1}(E_n)$ for all $n \geq 1$.\footnote{This implies that $M$ also has the structure of a Fredholm manifold, but modeled on $\E'$ \cite{Elworthy}.}
\end{dfn}

It follows from Theorem 4.3 of Gl\"ockner \cite{Glockner} that $M_\infty$ has the structure of a manifold based on the convenient vector space $\R^\infty = \indlim \R^n$, but we will not need this fact.

\begin{lem} Condition $(e)$ implies condition $(c)$ above. \end{lem}

\begin{proof} Since $M$ and each $M_n$ are separable Banach manifolds, they are ANRs, and so have the homotopy type of a CW-complex \cite{Milnor59,Palais}. Hence, the inductive limit $M_\infty = \indlim M_n$ has the homotopy type of a $CW$-complex by the Corollary on page 153 of \cite{Milnor63}. Condition (c) now follows by Whitehead's Theorem, since condition (d) shows that the induced maps on homotopy groups (which commute with direct limits) are isomorphisms. Indeed, the inverse to the induced map $(i_\infty)_* : \pi_k(M_\infty) \to \pi_k(M)$
is given by the following composition (using the open inclusions $U_n \into Z_n$):
$$\pi_k(M) = \pi_k(\cup_n U_n) = \indlim \pi_k(U_n) \to \indlim \pi_k(Z_n) \cong \indlim \pi_k(M_n) = \pi_k(M_\infty)$$ since each tubular neighborhood $Z_n$ strong deformation retracts onto $M_n$. \end{proof}

The following easy result, which is false for classical Fredholm filtrations,  shows the flexibility of these generalized $\Delta$-filtrations.

\begin{lem}\label{subseqlem} Any subsequence $\tilde{\F} = \{M_{n_k}\}_{k=1}^\infty$ of a (dense) normal $\Delta$-filtration is a (dense) normal $\tilde{\Delta}$-filtration where $\tilde{\Delta} = \{\delta(n_k)\}_{k=1}^\infty$ is the dimension sequence. Moreover, if $\F$ is Fredholm, then $\tilde{\F}$ is also Fredholm. \end{lem}

\begin{proof} Cofinality and normal bundles add and Lemma \ref{lem:flags} (i). \end{proof}

Now for some examples. 

\begin{exs} Let $\{E_n\}_1^\infty$ be a $\Delta$-flag in the Banach space $\E$. 
\begin{itemize}
\item[i.)] $\F = \{E_n\}_1^\infty$ is a dense normal $\Delta$-filtration  of $M = \E$ and is Fredholm via the identity map of $\E$.
\item[ii.)] If $U$ is an open subset of $\E$ such that $U \cap E_1 \neq \emptyset$ then $\F_U = \{U \cap E_n\}_1^\infty$ is a dense normal $\Delta$-filtration of $U$. The normal bundle of $U_n$ in $U$ is $\nu(U_n, U) = U_n \times E^{\infty -n}$. In fact, it is Fredholm via the inclusion map $U \into \E$.
\item[iii.)] Let $S(\E)$ denote the unit sphere of $\E$. Then, by Theorem 13.14 \cite{KrieglMichor}, $S(\E)$ is a smooth submanifold of $\E$ if and only if the norm on $\E$ is smooth away from the origin. For each $n \geq 1$, let $S^{\delta(n)-1} = S(\E) \cap E_n$ be the unit sphere in $E_n$. Then the sequence $\F_S = \{S^{\delta(n)-1}\}_1^\infty $ of finite-dimensional spheres is a dense normal $(\Delta-1)$-filtration of $S(\E)$. It is Fredholm via the inclusion map $S(\E) \into \E$. Note that inductive limit space $$M_\infty = \indlim S^{\delta(n) -1} \cong \indlim S^n = \S^\infty$$ is the infinite sphere.
\item[iv.)] If $M$ is a Fredholm manifold, then for some $k \geq 1$, if $\Delta_k = \{k, k+1, k+2, \dots \}$, a classical Fredholm filtration is a dense normal Fredholm $\Delta_k$-filtration and always exists due to Theorems 2.2 and 2.3 of \cite{Mkhr} if $\E$ has an $\N$-flag. See also Theorem 1.1 of \cite{Mkhr3}.
\item[v.)] Let $\F = \{M_n\}_1^\infty$ be a dense normal $\Delta$-filtration of $M$. Let $N = M \times \R^n$, then $\F_0 = \{M_n \times \{0\}\}_1^\infty$ induces a normal $\Delta$-filtration of $N$ that is not dense and not Fredholm.
\item[vi.)] Let $\F = \{M_n\}_1^\infty$ be a normal $\Delta$-filtration of $M$. Let $M' = M \times \R P^K$ where $\R P^K$ is real projective $K$-space for some $K > > 2$. Then $\F' = \{M'_n\}_1^\infty$ where$$M'_n =  \begin{cases} M_n \times \R P^n, 1 \leq n \leq K \\
M_n \times \R P^K, n > K \end{cases}$$ is a $\Delta'$-filtration of $M'$ that is neither normal nor Fredholm and is not dense, if $\F$ is not dense.
\item[vii.)] Let $\E$ be a real separable infinite-dimensional Hilbert space. Let $\P(\E)$ denote the projective space of $\E$, i.e., the quotient of the smooth (contractible) unit sphere $S(\E)$ by identifying antipodal points $x \sim -x$, which is a smooth Hilbert manifold. Using an $\N$-flag in $\E$ (e.g., from an orthonormal basis) we have canonical embeddings of real projective spaces
$$\RP^1 \subset \RP^2 \subset \RP^3 \subset \cdots \subset \RP^n \subset \cdots \subset  P(\E).$$
Thus, we have an inclusion of the infinite real projective space
$$\RP^\infty = \indlim \RP^n \into P(\E).$$
Using the fiber bundle $S^0 \to S(\E) \to P(\E)$  and Whitehead's Theorem (as above) it follows that the inclusion $\RP^\infty \into P(\E)$ is a (weak) homotopy equivalence. Thus, $\F = \{\RP^n \}_{n=1}^\infty$ is an $\N$-filtration of $P(\E)$ that is not normal and not Fredholm, due to the lack of trivial normal bundles, for any possible choice of Fredholm structure on $P(\E)$. A simple argument shows that $\RP^\infty = \S^\infty/\!\sim$ is dense in $P(\E) = S(\E)/\!\sim$ since $\S^\infty$ is dense in the unit sphere $S(\E)$.
\end{itemize}
\end{exs}

Combining example (iv) above with Lemma \ref{subseqlem}, gives the following existence result.  If $\E$ is a Hilbert space then the second statement follows by Theorem 3E \cite{EeElw2}.

\begin{prop} Let $M$ be a Fredholm manifold modelled on $\E$ with an $\N$-flag. Then for some $k \geq 1$ and any dimension sequence $\Delta = \{\delta(n)\}_1^\infty$ with $\delta(1) \geq k$, there is a Fredholm $\Delta$-filtration of $M$. In addition, if $\E$ is a Hilbert space, then we may assume that the tubular neighborhoods satisfy $V_n = U_n$ in condition (e).
\end{prop}

We now examine the functorial properties of these generalized $\Delta$-filtrations, which is not allowed in the classical theory due to the restrictive density and dimension requirements.

\begin{thm} If $\F = \{M_n\}_1^\infty$ is a (dense) normal $\Delta$-filtration of a Banach manifold $M$, then $\F_2 = \{ M_n\times M_n \}_{n=1}^\infty$ is a (dense) normal $2\Delta$-filtration of the pair groupoid $M\times M$.  Also, $(M \times M)_\infty \cong M_\infty \times M_\infty$ as topological spaces. If $\F$ is Fredholm, then so is $\F_2$.  \end{thm}

\begin{proof} Conditions (a), (b), and (d) are trivial. Condition (e) follows by using the tubular neighborhood $V_n \times V_n$ of $M_n \times M_n$ in $M \times M$ along with $U_n \times U_n$ as the associated open set. Since the $M_n$'s are locally compact Hausdorff, we have that $$(M \times M)_\infty = \indlim\  (M_n \times M_n) \cong M_\infty \times M_\infty$$ by Theorem 4.1 \cite{HHNH01}, i.e., the inductive limit and product topologies commute. If $M_\infty$ is dense in $M$, then it follows that $(M \times M)_\infty \cong M_\infty \times M_\infty$ is also dense in $M \times M$.  In the Fredholm case, if $f : M \to \E$ is the Fredholm map of index zero that induces $M_n = f^{-1}(E_n)$ relative to a $\Delta$-flag $\{E_n\}$ of $\E$ then the map $f \times f : M \times M \to \E \times \E$ is also a Fredholm map of index zero that is transversal to the $2\Delta$-flag  $\{E_n \times E_n\}_{n=1}^\infty$ in $\E \times \E$ and gives $(f \times f)^{-1}(E_n \times E_n) = M_n \times M_n$ and so the result follows. \end{proof}

Note that there is another filtration of $M \times M$:
$$M_1 \times M_1 \subset M_1 \times M_2 \subset M_2 \times M_2 \subset M_2 \times M_3 \subset M_3 \times M_3 \subset \cdots $$
but this is not a filtration of the pair groupoid $M \times M$ by pair (sub)groupoids. A slight generalization of the previous result easily follows.

\begin{thm} Let $M$ and $N$ be Banach manifolds (modelled on possibly different Banach spaces). If $\F_M = \{M_n\}_1^\infty$ is a (dense) normal $\Delta$-filtration of $M$ and $\F_N = \{N_n\}_{n=1}^\infty$ is a (dense) normal $\Delta'$-filtration of $N$, then $\F = \{ M_n\times N_n \}_{n=1}^\infty$ is a (dense) normal  $(\Delta + \Delta')$-filtration of $M\times N$. Also, $(M \times N)_\infty \cong M_\infty \times N_\infty$ as topological spaces. If $\F_M$ and $\F_N$ are  Fredholm, then so is $\F$. \end{thm}

In the classical Fredholm manifold literature, the tangent bundle $TM$ of a Fredholm manifold $M$ is stated to have an induced Fredholm manifold structure, but Fredholm filtrations of $TM$ are never discussed since the requirement that the dimension jumps by one is too restrictive. Note that the density requirement is dropped here since it is not clear that it would be true, but we have no counterexample in the non-Fredholm case. 

\begin{thm}\label{thm:TMfilt} If $\F = \{M_n\}_1^\infty$ normal $\Delta$-filtration of a Banach manifold $M$, then $T\F = \{ TM_n \}_{n=1}^\infty$ is a normal $2\Delta$-filtration of the tangent bundle $TM$. If $\F$ is Fredholm, then so is $T\F$. \end{thm}

\begin{proof} Condition (a) holds since $\dim(TM_n) = 2 \dim(M_n)$. The inclusions $M_n \subset M_{n+1}$ induce inclusions $TM_n \subset TM_{n+1}$ as manifolds so condition (b) holds. Condition (d) holds since if $p : NM_n \to M$ is a trivial normal bundle for $M_n$ in $M_{n+1}$ (or in $M$) then we have that $Tp: TNM_n \to TM_n$ is a trivial normal bundle for $TM_n$ in $TM_{n+1}$ (or $TM$), i.e., the tangent bundle and normal bundle functors commute \cite{BLM}. Condition (e) holds since $TV_n$ is a tubular neighborhood of $TM_n$ in $TM$ and the open sets $TU_n \subseteq TV_n$ satisfy $TU_n \subseteq TU_{n+1}$ and $TM = \cup_{n=1}^\infty TU_n$. If  $f : M \to \E$ is the Fredholm map of index zero that induces $M_n = f^{-1}(E_n)$ relative to a $\Delta$-flag $\{E_n\}$ of $\E$, then the differential $Df : TM \to T\E \cong \E \times \E$ is also a Fredholm map of index zero by Proposition \ref{p:dncmapfredholm} that is transversal to the $2\Delta$-flag $\{TE_n \cong E_n \times E_n\}_1^\infty$ by Corollary \ref{cor:trasveralitytangent} and satisfies $$(Df)^{-1}(E_n \times E_n) \cong (Df)^{-1}(TE_n) = Tf^{-1}(E_n)= TM_n$$ (which also follows by the Transversal Mapping Theorem 3.5.12 \cite{AMR88}.) \end{proof}
 
\begin{note} For Banach manifolds $M$, there are issues with the cotangent bundle $T^*M$, i.e., the dual vector bundle of the tangent bundle $TM$, even when $M = \E$. The first problem is that even if $\E$ is separable, the dual space $\E^*$ may not be separable (e.g., when  $\E = \ell^1$ we have $\E^* = \ell^\infty$.) Thus, $\E^*$ has no $\Delta$-flags. Moreover, even if $\E^*$ is separable, it does not follow that a $\Delta$-flag $\{E_n\}_1^\infty$ induces a $\Delta$-flag $\{E_n^*\}$ of the dual space $\E^*$. If one restricts to Hilbert Riemannian manifolds, then a cotangent bundle version of the previous theorem holds by taking vector bundle isomorphisms with the tangent bundle.
\end{note}
    
Now, we come to the main result, that led to the definition of generalized filtrations.

\begin{thm} If $\F = \{M_n\}_1^\infty$ is a normal $\Delta$-filtration of a Banach manifold $M$, then the sequence of tangent groupoids $\TG{\F} = \{ \TG{M_n} \}_{n=1}^\infty$ is a normal $(2\Delta+1)$-filtration of the tangent groupoid $\TG{M}$ of $M$. Moreover, if $\F$ is Fredholm, then so is $\TG{\F}$. \end{thm}

\begin{proof} Condition (a) holds since $\dim(\TG{M_n}) = 2 \dim(M_n) +1$. Condition (b) holds since the inclusions $M_n \subset M_{n+1}$ induce inclusions $\TG{M_n} \subset \TG{M_{n+1}}$ as manifolds. Condition (d) holds since if $p : NM_n \to M$ is a trivial normal bundle for $M_n$ in $M_{n+1}$ (or in $M$) then we have that $\TG{p}: \TG{NM_n} \to \TG{M_n}$ is a trivial normal bundle for $\TG{M_n}$ in $\TG{M_{n+1}}$ (or $\TG{M}$), by Proposition \ref{prop:trivialbdl}. Condition (e) holds since $\TG{V_n}$ is a tubular neighborhood of $\TG{M_n}$ in $\TG{M}$ and the open sets $\TG{U_n} \subseteq \TG{V_n}$ satisfy $\TG{U_n} \subseteq \TG{U_{n+1}}$ and $\TG{M} = \cup_{n=1}^\infty \TG{U_n}$. If  $f : M \to \E$ is a Fredholm map of index zero that induces $M_n = f^{-1}(E_n)$ relative to a $\Delta$-flag $\{E_n\}_1^\infty$ of $\E$, then the functorially induced map $\DG{f} : \TG{M} \to \TG{\E} \cong \E \times \E \times \R$ is also a Fredholm map of index zero by Proposition \ref{p:dncmapfredholm} that is transversal to the ($2\Delta$+1)-flag $\{\TG{E_n} \cong  E_n \times E_n \times\R\}_1^\infty$ by Corollary \ref{cor:transversaltangentgrpoid} and satisfies $$(\DG{f})^{-1}(E_n \times E_n \times \R) \cong (\DG{f})^{-1}(\TG{E_n}) = \TG{f^{-1}(E_n)} = \TG{M_n}$$ and the result now follows. \end{proof}

We finish this section with a few more easy results.

\begin{lem} If $f : M \to N$ and $g : N \to P$ are smooth maps of Banach manifolds such that $g$ is transversal to a submanifold $Z$ of $P$, then $f$ is transversal to $g^{-1}(Z)$ if and only the composition $g \circ f$ is transversal to $Z$.
\end{lem}

The following lemma is Theorem 3.3 of Smale \cite{SmaleSard}.

\begin{lem} Let $f : N \to M$ be a smooth Fredholm map between Banach manifolds with index $p \geq 0$ that is transverse to a finite-dimensional submanifold $M_n$ of $M$.  Then $N_n = g^{-1}(M_n)$ is a submanifold of $N$ of dimension $\dim(N_n) = p + \dim(M_n).$
\end{lem}

Combining these two lemmas shows how to induce generalized Fredholm filtrations with dimension shifts via Fredholm maps with higher indices into $M$. Compare the discussion after Theorem 1.1 \cite{Mkhr3}. 

\begin{thm} Let $\F = \{M_n\}_1^\infty$ be a Fredholm $\Delta$-filtration of the Fredholm manifold $M$. If $N$ is a Fredholm manifold and $g : N \to M$ is a Fredholm map of index $p \geq 0$ that is transverse to each $M_n$, then $N_n = g^{-1}(M_n)$ induces a Fredholm $(\Delta+p)$-filtration $\F_g = \{N_n\}_1^\infty$ of $N$. \end{thm}

Note we assume that $N_1 = g^{-1}(M_1) \neq \emptyset$ which we can do by Lemma \ref{subseqlem}. The following is then easy to prove.

\begin{prop} Let $p : \widetilde{M} \to M$ be a smooth covering of the Banach manifold $M$ and let \F = $\{M_n\}_1^\infty$ be a (dense) normal $\Delta$-filtration of $M$. For each $n \geq 1$, let $\widetilde{M}_n = p^{-1}(M_n)$. Then $\widetilde{\F} = \{\widetilde{M_n}\}_1^\infty$ is a (dense) normal $\Delta$-filtration of the Banach manifold $\widetilde{M}$ by covering spaces $p_n = p|_{M_n} : \widetilde{M_n} \to M_n$. Moreover, if $\F$ is Fredholm, then so is $\widetilde{\F}$.
\end{prop}

\bibliographystyle{elsarticle-num}
\bibliography{Revised.master}

\end{document}